\def\version{July 10, 2019}

\documentclass[11pt]{article}
\usepackage{amsmath,enumerate, amsfonts, amssymb, amsthm, bbm}
\usepackage{xcolor}

\newcommand{\ch}[1]{#1}

\def\UseSection{
        \numberwithin{equation}{section}
    \theoremstyle{plain}
        \newtheorem{theorem}    {Theorem}[section]
        \DefineTheorems 
}

\def\DefineTheorems{
    
    \newtheorem{lemma}      [theorem] {Lemma}
    
    \newtheorem{prop}       [theorem] {Proposition}
    
    \newtheorem{cor}        [theorem] {Corollary}

    \theoremstyle{definition}
    \newtheorem{defn}       [theorem] {Definition}

    \theoremstyle{definition}

}

\newcommand{\bt}   {\begin{theorem}}
\newcommand{\et}   {\end  {theorem}}
\newcommand{\bl}   {\begin{lemma}}
\newcommand{\el}   {\end  {lemma}}
\newcommand{\bp}   {\begin{prop}}
\newcommand{\ep}   {\end  {prop}}
\newcommand{\bc}   {\begin{cor}}
\newcommand{\ec}   {\end  {cor}}
\newcommand{\bd}   {\begin{defn}}
\newcommand{\ed}   {\end  {defn}}

\newcommand{\ba}   {\begin{array}}
\newcommand{\ea}   {\end  {array}}
\newcommand{\be}   {\begin{enumerate}}
\newcommand{\ee}   {\end  {enumerate}}
\newcommand{\bi}   {\begin{itemize}}
\newcommand{\ei}   {\end  {itemize}}

\def\eq#1\en{\begin{equation}#1\end{equation}}
\def\eqsplit#1\ensplit{
    \begin{equation}\begin{split}#1\end{split}\end{equation}
    }
\def\eqalign#1\enalign{
    \begin{align}#1\end{align}
    }
\def\eqmul#1\enmul{
    \begin{multline}#1\end{multline}
    }
\newcommand{\eqarrstar} {\begin{eqnarray*}}
\newcommand{\enarrstar} {\end{eqnarray*}}
\newcommand{\eqarray}   {\begin{eqnarray}}
\newcommand{\enarray}   {\end{eqnarray}}
\newcommand{\nnb}   {\nonumber \\}

\newcommand{\lbeq}[1]  {\label{e:#1}}
\newcommand{\refeq}[1] {\eqref{e:#1}}    

\makeatletter
\newcommand{\labelcounter}[2]{{%
    \stepcounter{#1}
    \protected@write\@auxout{}%
    {\string\newlabel{#2}{{\csname the#1\endcsname}{\thepage}}}%
    {\ref{#2}}
    }}
\makeatother
\newcommand{\sss}   { \scriptscriptstyle }

\newcommand{\Ebold} {{\mathbb E}}

\newcommand{\Pbold} {{\mathbb P}}

\newcommand{\Zbold} {{\mathbb Z}}

\newcommand{\Ccal}   {\mathcal{C}}

\newcommand{\Zd}    {{ {\Zbold}^d }}

\newcommand{\spose}[1] {{\hbox to 0pt{#1\hss}} }
\newcommand{\ltapprox} {\mathrel{\spose{\lower 3pt\hbox{$\mathchar"218$}}
 \raise 2.0pt\hbox{$\mathchar"13C$}}}
\newcommand{\gtapprox} {\mathrel{\spose{\lower 3pt\hbox{$\mathchar"218$}}
 \raise 2.0pt\hbox{$\mathchar"13E$}}}

\newcommand{\expec}[1]  {\left \langle #1 \right \rangle}

\UseSection  
\newcommand{\C}     {\mathcal {C}}

\newcommand{\N}     {\mathbb{N}}
\renewcommand{\P}   {\mathbb{P}}

\newcommand{\E}     {\mathbb{E}}
\newcommand{\Z}     {\mathbb{Z}}

\newcommand\1{\mathbbm{1}}

\newcommand{\Ci}{\Ccal_{\sss(i)}}
\newcommand{\Cone}{\Ccal_{\sss(1)}}
\newcommand{\Ctwo}{\Ccal_{\sss(2)}}
\newcommand{\diam}{{\rm diam}}

\newcommand{\sZ}{{\sss \mathbb{Z}}}
\newcommand{\sT}{{\sss \mathbb{T}}}
\newcommand{\torus}{\mathbb T_{r,d}}

\newcommand{\gr}{{\mathbb G}}

\newcommand{\V}{\mathbb V}

\newcommand{\ver}{{\mathcal V}}
\newcommand{\edges}{{\mathcal E}}
\newcommand{\graph}{{\mathcal G}}
\newcommand{\Cmax}{{\cal C}_{\rm max}}

\newcommand{\egr}\Exp
\newcommand{\cn}{\Omega}

\newcommand{\Cg}{\tilde{C}_{\sss\rm\chi}}
\newcommand{\cxi}{c_{\rm{\sss \xi}}}
\newcommand{\Cxi}{C_{\rm{\sss \xi}}}
\newcommand{\cxit}{c_{\rm{\sss \xi_2}}}
\newcommand{\Cxit}{C_{\rm{\sss \xi_2}}}

\newcommand{\Cwindow}{C_{\sss p_c}}
\newcommand{\cdelta}{c_{\rm{\sss \Ccal}}}

\newcommand{\Ctilde}{{\tilde C}}

\renewcommand{\expec}{{\mathbb E}}

\newcommand{\Var}{{{\rm Var}_p}}
\newcommand{\VarPc}{{{\rm Var}_{p_c(\Zd)}}}
\newcommand{\Exp}{{{\mathbb E}_p}}

\newcommand{\prob}{{\mathbb P}}

\newcommand{\conn}{\leftrightarrow}
\newcommand{\nc}{\conn{\hspace{-2.8ex} /}\hspace{1.8ex}}
\newcommand{\nn}{\nonumber}

\newcommand{\eqn}[1]{\begin{equation}#1\end{equation}}
\newcommand{\eqan}[1]{\eqalign #1 \enalign}

\newcommand{\indic}[1]{\1_{\{#1\}}}

\newcommand{\ptrans}{{\textbf{\rm  p}}}
\newcommand{\Tmix}{T_{{\rm mix}}}

\newcommand{\e}{{\mathrm e}}

\begin{document}

\thispagestyle{empty}
\vspace{1cm}
\centerline{\LARGE \bf Random graph asymptotics on high-dimensional tori}
\vspace{0.2cm}
\centerline{\LARGE \bf II. Volume, diameter and mixing time}
\vspace{1cm}

\centerline {{\sc Markus Heydenreich$^1$} and {\sc Remco van der Hofstad$^2$}}
\vspace{.6cm}

\centerline{\em $^1$Vrije Universiteit Amsterdam, Department of Mathematics,}
\centerline{\em De Boelelaan 1081a, 1081 HV Amsterdam, The~Netherlands}
\centerline{\tt MO.Heydenreich@few.vu.nl}
\vspace{.3cm}

\centerline{\em $^2$Eindhoven University of Technology,}
\centerline{\em Department of Mathematics and Computer Science,}
\centerline{\em P.O.~Box 513, 5600~MB Eindhoven, The~Netherlands}
\centerline{\tt r.w.v.d.hofstad@tue.nl}
\vspace{1cm}

\centerline{\small (\version)}
\vspace{.6cm}

\begin{quote}
  {\small {\bf Abstract:}}
  For critical bond-percolation on high-dimensional torus, this paper proves sharp lower bounds on the size of the largest cluster, removing a logarithmic correction in the lower bound in \cite{HeydeHofst07}.
  This improvement finally settles a conjecture by Aizenman \cite{Aizen97} about the role of boundary conditions in critical high-dimensional percolation, and it is a key step in deriving further properties of critical percolation on the torus.
  Indeed, a criterion of Nachmias and Peres \cite{NachmPeres08} implies appropriate bounds on diameter and mixing time of the largest clusters.
  We further prove that the volume bounds apply also to any finite number of the largest clusters.
  The main conclusion of the paper is that the behavior of critical percolation on the high-dimensional torus is the same as for critical Erd\H{o}s-R\'enyi random graphs.

In this updated version we incorporate an erratum to be published  in a forthcoming issue of \emph{Probab.\ Theory Relat.\ Fields}. This results in a modification of Theorem 1.2 as well as Proposition 3.1. 
\end{quote}

\vspace{0.5cm}


\vspace{1cm}
\hrule
\vspace{1em}

\section{Introduction}
\subsection{The model}
For bond percolation on a graph $\gr$ we make any edge (or `bond') \emph{occupied} with probability $p$, independently of each other, and otherwise leave it \emph{vacant}.
The connected components of the random subgraph of occupied edges are called \emph{clusters}.
For a vertex $v$ we denote by $\Ccal(v)$ the unique cluster containing $v$, and by $|\Ccal(v)|$ the number of vertices in that cluster.
For our purposes it is important to consider clusters as subgraphs (thus not only as a set of vertices).
Our main interest is bond percolation on high-dimensional tori, but our techniques are based on a comparison with $\Zd$ results.
We describe the $\Zd$-setting first.

\paragraph{Bond percolation on $\Zd$.}
For $\gr=\Zd$, we consider two sets of edges.
In the {\em nearest-neighbor model}, two vertices $x$ and $y$ are linked by an edge whenever $|x-y|=1$, whereas in the {\em spread-out model}, they are linked whenever $0<\|x-y\|_\infty\le L$.
Here, and throughout the paper, we write $\|\cdot\|_\infty$ for the supremum norm, and $|\cdot|$ for the Euclidean norm.
The integer parameter $L$ is typically chosen large.

The resulting product measure for percolation with parameter $p\in[0,1]$ is denoted by $\Pbold_{\sZ,p}$, and the corresponding expectation $\Ebold_{\sZ,p}$.
We write $\{0\conn x\}$ for the event that there exists a path of occupied edges from the origin $0$ to the lattice site $x$ (alternatively, $0$ and $x$ are in the same cluster), and define
    \eq
    \lbeq{tau-def}
    \tau_{\sZ,p}(x):=\Pbold_{\sZ,p}(0\conn x)
    \en
to be the \emph{two-point} function.
By
$$\chi_\sZ(p):=\sum_{x\in\Zd}\tau_{\sZ,p}(x)=\Ebold_{\sZ,p}|\Ccal(0)|$$
we denote the expected cluster size on $\Zd$.
The degree of the graph, which we denote by $\cn$, is $\cn=2d$ in the nearest-neighbor case and $\cn=(2L+1)^d-1$ in the spread-out case.

Percolation on $\Zd$ undergoes a phase transition as $p$ varies, and it is well known that there exists a critical value
    \eqn{
    \lbeq{pc-Zd}
        p_c(\Zd)
        =\inf\{p\colon \Pbold_{\sZ,p}(|\Ccal(0)|=\infty)>0\}
        =\sup\{p\colon \chi_\sZ(p)<\infty\},
    }
\ch{where the last equality is due to Aizenman and Barsky \cite{AizenBarsk87} and Menshikov \cite{Mensh86}.}

\paragraph{Bond percolation on the torus.}
By $\torus$ we denote a graph with vertex set
$\{-\lfloor r/2\rfloor,\dots,\lceil r/2\rceil-1\}^d$ and two related sets of edges:
\begin{enumerate}
\item The nearest-neighbor torus:
an edge joins vertices that differ by $1$ (modulo $r$)
in exactly one component.
For $d$ fixed and $r$ large, this is a periodic approximation to $\Z^d$.
Here $\Omega =2d$ for $r \geq 3$.
We study the limit in which $r\rightarrow \infty$ with $d>6$ fixed, but large.

\item The spread-out torus:
an edge joins vertices $x=(x_1,\ldots,x_d)$ and $y=(y_1,\ldots,y_d)$
if $0<\max_{i=1,\ldots,d}|x_i-y_i|_r \leq L$ (with $|\cdot|_r$ the metric
on $\Zbold_r$).  We study the limit $r \to \infty$,
with $d > 6$ fixed and $L$ large (depending on $d$)
and fixed.  This gives a periodic approximation to
range-$L$ percolation on $\Z^d$. Here $\Omega = (2L+1)^d-1$
provided that $r\geq 2L+1$, which we will always assume.
\end{enumerate}
We write $V=r^d$ for the number of vertices in the torus.
We consider bond percolation on these tori with edge occupation probability $p$ and write $\Pbold_{\sT,p}$ and
$\Ebold_{\sT,p}$ for the product measure and corresponding expectation, respectively.
We use notation analogously to $\Zd$-quantities, e.g.\
$$\chi_\sT(p):=\sum_{x\in\torus}\Pbold_{\sT,p}(0\conn x)=\Ebold_{\sT,p}|\Ccal(0)|$$
for the expected cluster size on the torus.

\ch{\paragraph{Mean-field behavior in high dimensions.} In the past decades, there has
been substantial progress in the understanding of percolation in high-dimensions
(see e.g. \cite{AizenNewma84, BarskAizen91, Hara90,Hara08,HaraHofstSlade03,HaraSlade90a, HaraSlade00a, HaraSlade00b, Slade06} for detailed results on high-dimensional
percolation), and the results show that percolation
on high-dimensional infinite lattices is similar to percolation on infinite trees
(see e.g., \cite[Section 10.1]{Grimm99} for a discussion of percolation on a tree).
Thus, informally speaking, the mean-field model for percolation on $\Z^d$
is percolation on the tree.

More recently, the question has been addressed what the mean-field model is of
percolation on \emph{finite} subsets of $\Z^d$, such as the torus.
Aizenman \cite{Aizen97} conjectured that critical percolation on high-dimensional
tori behaves similarly to critical Erd\H{o}s-R\'enyi random graphs, thus suggesting
that the mean-field model for percolation on a torus is the
Erd\H{o}s-R\'enyi random graph. In the past years, substantial progress was
made in this direction, see in particular \cite{BorgsChayeHofstSladeSpenc05a, BorgsChayeHofstSladeSpenc05b,HeydeHofst07}. In this paper, we bring this discussion
to the next level, by showing that large critical clusters on various
high-dimensional tori} share many features of the Erd\H{o}s-R\'enyi random graph.

\subsection{Random graph asymptotics on high-dimensional tori}
\label{sec-high-d}
We investigate the size of the maximal cluster on the torus $\torus$, i.e.,
    \eq
    |\Cmax|:=\max_{x\in\torus}|\Ccal(x)|,
    \en
at the critical percolation threshold $p_c(\Zd)$.
We start by improving the asymptotics of the largest connected component as proved in
\cite{HeydeHofst07}:

\begin{theorem}[Random graph asymptotics of the largest cluster size]
\label{thm-1}
Fix $d>6$ and $L$ sufficiently large in the spread-out case, or $d$ sufficiently large for nearest-neighbor
percolation.
Then there exists a constant $b>0$, such that for all $\omega\ge1$ and all $r\ge1$,
    \eq
    \lbeq{Cmaxbd}
    \Pbold_{\sT, p_c(\Z^d)}\Big(\omega^{-1}V^{2/3}
    \leq |\Cmax|\leq \omega V^{2/3}\Big)
    \geq 1-\frac{b}{\omega}.
    \en
The constant $b$ can be chosen equal to $b_6$ in \cite[Theorem 1.3]{BorgsChayeHofstSladeSpenc05a}.
Furthermore, there are positive constants $c_1$ and $c_2$ such that
    \eq
    \lbeq{Cmaxbd2}
    \Pbold_{\sT, p_c(\Z^d)}\Big(|\Cmax|> \omega V^{2/3}\Big)
    \le \frac{c_1}{\omega^{3/2}}\;{\rm e}^{-c_2\,\omega}.
    \en
\end{theorem}
We recall that $r$ is present in \refeq{Cmaxbd} in two ways: We consider the percolation measure on $\torus$, and $V=r^d$ is the volume of the torus.
The upper bound \ch{in \refeq{Cmaxbd} in Theorem \ref{thm-1}} is already
proved in \cite[Theorem 1.1]{HeydeHofst07},
whereas the lower bound in \cite[Theorem 1.1]{HeydeHofst07} contains a
logarithmic correction, which we remove
here by a more careful analysis.

We next extend the above result to the other large clusters. For this, we write
$\Ci$ for the $i^{\rm th}$ largest cluster for percolation on $\torus$, so that
$\Cone=\Cmax$ and $|\Ctwo|\leq |\Cone|$ is the size of the second largest component; etc.

\begin{theorem}[Random graph asymptotics of the ordered cluster sizes]
\label{thm-2}
Fix $d>6$ and $L$ sufficiently large in the spread-out case, or $d$ sufficiently large for nearest-neighbor
percolation.
For every $m=1,2,\dots$ there exist constants $b_1, \dots, b_m>0$, such that for all $\omega\ge1$, $r\ge1$,
and all $i= 1,\dots,m$,
    \eq
    \lbeq{Cibd}
    \Pbold_{\sT, p_c(\Z^d)}\Big(\omega^{-1}V^{2/3}
    \leq |\Ci|\leq \omega V^{2/3}\Big)
    \geq 1-\frac{b_i}{\omega}.
    \en
Consequently, the expected cluster sizes satisfy $\E_{\sT,p_c(\Zd)}|\Ci|\ge b_i'\,V^{2/3}$ for certain constants $b_i'>0$.
\end{theorem}

Nachmias and Peres \cite{NachmPeres08} proved a very handy criterion establishing bounds on \emph{diameter} and \emph{mixing time of lazy simple random walk} of the large critical clusters for random graphs obeying \eqref{e:Cmaxbd}/\eqref{e:Cibd}.
The following corollary states the consequences of the criterion for the high-dimensional torus.
To this end, we call a \emph{lazy simple random walk} on a finite graph $\graph=(\ver,\edges)$ a Markov chain on the vertices $\ver$ with transition probabilities
\begin{equation}\label{eq:defLazyRW}
    \ptrans(x,y)=
    \begin{cases}
    1/2\qquad&\text{if $x=y$;}\\
    \frac1{2\,\operatorname{deg}(x)}\qquad&\text{if $(x,y)\in\edges$;}\\
    0\qquad&\text{otherwise,}\\
    \end{cases}
\end{equation}
where $\operatorname{deg}(x)$ denotes the degree of a vertex $x\in\ver$.
The stationary distribution of this Markov chain $\pi$ is given by $\pi(x)=\operatorname{deg}(x)/(2|\edges|)$.
The \emph{mixing time} of lazy simple random walk is defined as
\begin{equation}\label{eq:defMixingTime}
    \Tmix(G)=\min\big\{n\colon\|\ptrans^n(x,\cdot)-\pi(\cdot)\|_{\rm \sss TV}\le 1/4 \text{ for all $x\in V$}\big\},
\end{equation}
with $\ptrans^n$ being the distribution after $n$ steps (i.e., the $n$-fold convolution of $\ptrans$), and $\|\cdot\|_{\rm \sss TV}$ denoting the total variation distance.
We write $\diam(\Ccal)$ for the diameter of the
cluster $\Ccal$.

\begin{cor}[{Diameter and mixing time of large critical clusters, \cite{NachmPeres08}}]
\label{thm-3}
Fix $d>6$ and $L$ sufficiently large in the spread-out case, or $d$ sufficiently large for nearest-neighbor
percolation.
Then, for every $m=1,2,\dots$, there exist constants $c_1, \dots, c_m>0$, such that for all $\omega\ge1$, $r\ge1$,
and all $i= 1,\dots,m$,
    \eqan{
    \lbeq{diamCibd}
    \Pbold_{\sT, p_c(\Z^d)}\Big(\omega^{-1}V^{1/3}
    \leq \diam(\Ci)\leq \omega V^{1/3}\Big)
    &\geq 1-\frac{c_i}{\omega^{1/3}},\\
    \label{eq:Tmix}
    \Pbold_{\sT, p_c(\Z^d)}\Big(\omega^{-1}V\leq\Tmix(\Ci)\leq \omega V\Big)
    &\geq 1-\frac{c_i}{\omega^{1/34}}.
    }
\end{cor}

\subsection{Discussion and open problems}
\label{sec-discussion}

Here, and throughout the paper, we make use of the following notation:
we write $f(x)=O(g(x))$ for functions $f,g\ge0$ and $x$ converging to some limit, if there exists a constant $C>0$ such that $f(x)\le Cg(x)$ in the limit, and $f(x)=o(g(x))$ if $g(x)\neq O(f(x))$.
Furthermore, we write  $f=\Theta(g)$ if $f=O(g)$ and $g=O(f)$.

The asymptotics of $|\Cmax|$ in Theorem \ref{thm-1} is an improvement of our earlier result in \cite{HeydeHofst07}, which itself relies in an essential way on the work of Borgs et al.\ \cite{BorgsChayeHofstSladeSpenc05a,BorgsChayeHofstSladeSpenc05b}.
The contribution of the present paper is the removal of the logarithmic correction in the lower bound of \cite[(1.5)]{HeydeHofst07}, and this improvement is crucial for our further results, as we discuss
in more detail now. We give an easy proof that the largest $m$ components obey the same volume
asymptotic as the largest connected component, using only Theorem \ref{thm-1} and estimates
on the moments of the random variable
    \eqn{
    \lbeq{Zgeq-def}
    Z_{\sss \geq k}=\#\{v\in\torus: |\Ccal(v)|\ge k\}
    }
derived in \cite{BorgsChayeHofstSladeSpenc05a,BorgsChayeHofstSladeSpenc05b}.
Given these earlier results, our proofs are remarkably simple and robust, and they can be expected to apply in various different settings.
Thus, while our results substantially improve our understanding of the critical nature of percolation on high-dimensional tori, the proofs given here are surprisingly simple.

\paragraph{Random graph asymptotics at criticality.}
Our results show that the largest percolation clusters on the high-dimensional torus behave as they do on the Erd\H{o}s-R\'enyi random graph; this can be seen as the take-home message of this paper.
Aldous \cite{Aldou97} proved that, for Erd\H{o}s-R\'enyi random graphs, the vector
$$V^{-2/3}\big(|\Cone|,|\Ctwo|,\dots,|\C_{(m)}|\big)$$
converges in distribution, as $V\to\infty$, to a random vector
    $(|\gamma_1|,\dots,|\gamma_m|)$,
where $|\gamma_j|$ are the excursion lengths (in decreasing order) of reflected Brownian motion.
Nachmias and Peres \cite[Thm.\ 5]{NachmPeres09b} prove the same limit
(apart from a multiplication with an explicit constant)
for random $d$-regular graphs (for which the critical value equals $(d-1)^{-1}$).
In light of our Theorems \ref{thm-1}--\ref{thm-2}, we conjecture that \ch{the same limit,
multiplied by an appropriate constant as in \cite[Thm.\ 5]{NachmPeres09b},
arises} for the ordered largest critical components for percolation
on high-dimensional tori.

\paragraph{The role of boundary conditions.}
The combined results of Aizenman \cite{Aizen97} and Hara et al.\ \cite{Hara08,HaraHofstSlade03} show that a box of width $r$ under \emph{bulk} boundary conditions in high dimension
satisfies $|\Cmax|\approx r^4$, which is much smaller than $V^{2/3}$.
This immediately implies an upper bound on $|\Cmax|$ under \emph{free} boundary conditions.
Aizenman \cite{Aizen97} conjectures that, under periodic boundary conditions, $|\Cmax|\approx V^{2/3}$.
This conjecture was proven in \cite{HeydeHofst07} with a logarithmic correction in the lower bound.
The present paper (improving the lower bound) is the ultimate confirmation of the conjecture in \cite{Aizen97}.

\paragraph{The critical probability for percolation on the torus.}
An alternative definition for the critical percolation threshold on a general high-dimensional torus, denoted by $p_c(\torus)$, was given in \cite[(1.7)]{BorgsChayeHofstSladeSpenc05a} as the solution to
    \eq\lbeq{defPcT}
    \chi_{\sT}(p_c(\torus))=\lambda V^{1/3},
    \en
where $\lambda$ is a sufficiently small constant.
The definition of the critical value in \refeq{defPcT} appears somewhat indirect, but the big advantage is that this definition exists for any torus (including $d$-cube, Hamming cube, complete graph), even if an externally defined critical value (such as $p_c(\Zd)$ as in \refeq{pc-Zd}) does not exist.
It is a major result of Borgs et al.\ \cite{BorgsChayeHofstSladeSpenc05a,BorgsChayeHofstSladeSpenc05b} that Theorem \ref{thm-1} holds with $p_c(\Zd)$ replaced by $p_c(\torus)$ for the following tori:
\begin{enumerate}
\item the $d$-cube ${\mathbb T}_{2,d}$ as $d\rightarrow \infty$,
\item the complete graph (Hamming torus with $d=1$ and $r \to \infty$),
\item
nearest-neighbor percolation on $\torus$ with $d \geq 7$
and $r^d \to \infty$ in any fashion, including $d$ fixed and $r \to \infty$,
$r$ fixed and $d \to \infty$, or $r,d \to \infty$ simultaneously,
\item periodic approximations to range-$L$ percolation on $\Z^d$ for
fixed $d \geq 7$ and fixed large $L$.
\end{enumerate}
Remarkably, our results in Theorem \ref{thm-2} and Corollary \ref{thm-3} hold also for all of the above listed tori when $p_c(\Zd)$ is replaced by $p_c(\torus)$.
One way of formulating Theorem \ref{thm-1} is to say that $p_c(\torus)$ and
$p_c(\Z^d)$, under the assumptions of Theorem \ref{thm-1}, are asymptotically equivalent.


\section{Proof of Theorem \ref{thm-1}}
\label{sec-thm-1}
The following relation between the two critical values $p_c(\Zd)$ (which is `inherited' from the infinite lattice) and $p_c(\torus)$ (as defined in \refeq{defPcT}) is crucial for our proof.
\begin{theorem}[The $\Zd$ critical value is inside the $\torus$ critical window]
\label{thm-ImprovedBoundsCmax}
Fix $d>6$ and $L$ sufficiently large in the spread-out case, or $d$ sufficiently large for nearest-neighbor
percolation.
Then there exists $\Cwindow>0$ such that $p_c(\Zd)$ and $p_c(\torus)$ satisfy
\begin{equation}\label{eq:pc-bound}
    \left|p_c(\Zd)-p_c(\torus)\right|\le \Cwindow \,V^{-1/3}.
\end{equation}
\end{theorem}

\noindent
In other words, $p_c(\Zd)$ lies in a critical window of order $V^{-1/3}$ around $p_c(\torus)$.
By the work of Borgs, Chayes, van der Hofstad, Slade and Spencer
\cite{BorgsChayeHofstSladeSpenc05a,BorgsChayeHofstSladeSpenc05b},
Theorem \ref{thm-ImprovedBoundsCmax} has immediate consequences for the size of the largest cluster, and various other quantities:
\begin{cor}[Borgs et al.\ \cite{BorgsChayeHofstSladeSpenc05a,BorgsChayeHofstSladeSpenc05b}]
\label{cor:BCHSS-bounds}
Under the conditions of Theorem \ref{thm-ImprovedBoundsCmax}, there exist constants $b,C>0$, such that for all $\omega\ge C$,
    \eq
    \label{eq:CmaxBound}
    \Pbold_{\sT, p_c(\Z^d)}\Big(\omega^{-1}V^{2/3}
    \leq |\Cmax|\leq \omega V^{2/3}\Big)
    \geq 1-\frac b {\omega}.
    \en
Furthermore,
\begin{equation}
    c\,V^{2/3}\le \E_{\sT, p_c(\Z^d)}\big(|\Cmax|\big)\le C\,V^{2/3}
\quad \text{and} \quad
\label{eq:ClusterExpect}
c_{\sss\chi}V^{1/3}\le\E_{\sT, p_c(\Z^d)}\big(|\C|\big)\le C_{\sss\chi}V^{1/3}
\end{equation}
for some $c,C, c_{\sss \chi}, C_{\sss \chi}>0$.
Finally, there are positive constants $b_{\sss\Ccal}$, $\cdelta$ 
such that for $k\le b_{\sss\Ccal}V^{2/3}$,
\begin{equation}\label{eq:DeltaEqualsTwo}
    \frac\cdelta{\sqrt{k}}\le \P_{\sT, p_c(\Z^d)}\big(|\C|\ge k\big) \le \frac6{\sqrt{k}}.
\end{equation}
All of these statements hold uniformly as $r\to\infty$.
\end{cor}
The reader may verify that Corollary \ref{cor:BCHSS-bounds} indeed follows from Theorem \ref{thm-ImprovedBoundsCmax} by using
\cite[Thm.\ 1.3]{BorgsChayeHofstSladeSpenc05a} in conjunction with
\cite[Prop.\ 1.2 and Thm.\ 1.3]{BorgsChayeHofstSladeSpenc05b}.
Note that \eqref{eq:CmaxBound} in particular proves \refeq{Cmaxbd} in Theorem \ref{thm-1}.
The explicit value of the constant in the upper bound of \eqref{eq:DeltaEqualsTwo} is not mentioned explicitly in \cite{BorgsChayeHofstSladeSpenc05a}, but an inspection of the proof of the upper bound in \cite[Theorem 1.3]{BorgsChayeHofstSladeSpenc05a} shows that $6/\sqrt k$ suffices. Indeed, by \cite[Proposition 2.1]{HeydeHofst07}, $\P_{\sT, p_c(\Zd)}(|\C|\geq k)\leq \P_{\sZ, p_c(\Zd)}(|\C|\geq k)$. Further, by \cite[(9.2.6)]{HeyHof17}, $\P_{\sZ, p_c(\Zd)}(|\C|\geq k)\leq \frac{\e}{\e-1}M(p_c(\Zd),1/k)$, while \cite[Lemma 9.3]{HeyHof17} proves that $M(p_c(\Zd), \gamma)\leq \sqrt{12\gamma}$.

We explicitly keep track of the origin of constants by adding an appropriate subscript.
For first time reading the reader might wish to ignore these subscripts.

We are now turning towards the proof of Theorem \ref{thm-ImprovedBoundsCmax}.
To this end, we need the following lemma:
\begin{lemma}\label{lemma:SumTauBound}
    For percolation on $\Zd$ with $p=p_c(\Z^d)-{K}{\cn}^{-1}V^{-1/3}$, there exists a positive constant $\Ctilde$ (depending on $d$ and $K$, but not on $V$), such that
    \begin{equation}\label{eq:chiBd}
        \sum_{\substack{u,v\in\Zd, u\neq v \\ u-v\in r\Zd}}\tau_p(u)\,\tau_p(v)\le \Ctilde\,V^{-1/3}.
    \end{equation}
\end{lemma}
The lemma makes use of a number of results on high-dimensional percolation on $\Zd$, to be summarized in the following theorem.
\begin{theorem}[$\Zd$-percolation in high dimension \cite{Hara90,Hara08,HaraHofstSlade03,HaraSlade90a}.]
\label{thm:ZdResults}
Under the conditions in Theorem \ref{thm-1},  there exist constants $c_{\sss \tau}, C_{\sss \tau}, \cxi, \Cxi, \cxit, \Cxit>0$ such that
    \eq
    \lbeq{powerbd}
    \frac{c_{\sss \tau}}{(|x|+1)^{d-2}} \leq \tau_{\sZ,p_c(\Zd)}(x) \leq \frac{C_{\sss \tau}}{(|x|+1)^{d-2}}.
    \en
Furthermore, for any $p<p_c(\Zd)$,
    \eq
    \lbeq{expbd}
    \tau_{\sZ,p}(x) \leq {\rm e}^{-\frac{\|x\|_\infty}{\xi(p)}},
    \en
where the \emph{correlation length} $\xi(p)$ is defined by
    \eq
    \lbeq{defXi}
    \xi(p)^{-1}=-\lim_{n\to\infty}\frac1n\log\Pbold_{\sZ,p}\big((0,\dots,0)\conn(n,0,\dots,0)\big),
    \en
and satisfies
    \eq
    \lbeq{xibd}
    \cxi\left(p_c(\Z^d)-p\right)^{-1/2}\le\xi(p)\le\Cxi\left(p_c(\Z^d)-p\right)^{-1/2}\qquad\text{as $p\nearrow p_c(\Z^d)$}.
    \en
For the \emph{mean-square displacement}
    \begin{equation}\label{eq:CorrelLengthOrderTwo}
    \xi_2(p):=\left(\frac{\sum_{v\in\Zd}|v|^2\tau_{\sZ,p}(v)}{\sum_{v\in\Zd}\tau_{\sZ,p}(v)}\right)^{1/2},
    \end{equation}
we have
    \eq
    \lbeq{xi2bd}
    \cxit\left(p_c(\Z^d)-p\right)^{-1/2}\le\xi_2(p)\le\Cxit\left(p_c(\Z^d)-p\right)^{-1/2}\qquad\text{as $p\nearrow p_c(\Z^d)$}.
    \en
Finally, there exists a positive constant $\Cg$, such that the \emph{expected cluster size} $\chi_\sZ(p)$ obeys
\eq\lbeq{defCg}
\frac1{\cn\left({p_c(\Z^d)-p}\right)}\le\chi_\sZ(p)\le\frac\Cg{\cn\left({p_c(\Z^d)-p}\right)}\qquad\text{as $p\nearrow p_c(\Z^d)$}.
\en
\end{theorem}
Some of these bounds express that certain critical exponents exist and take on their mean-field value.
For example, \refeq{powerbd} means that the $\eta=0$, and similarly \refeq{defCg} can be rephrased as $\gamma=1$.
The power-law bound \refeq{powerbd} is due to Hara \cite{Hara08} for the nearest-neighbor case, and to Hara, van der Hofstad and Slade \cite{HaraHofstSlade03} for the spread-out case.
For the exponential bound \refeq{expbd}, see e.g.\ Grimmett \cite[Prop.\ 6.47]{Grimm99}.
Hara \cite{Hara90} proves the bound \refeq{xibd}, and Hara and Slade \cite{HaraSlade90a} prove \refeq{xi2bd} and \refeq{defCg} (the latter in conjunction with Aizenman and Newman \cite{AizenNewma84}).
The proof of all of the above results uses the lace expansion.

\proof[Proof of Lemma \ref{lemma:SumTauBound}.]
We split the sum on the left-hand side of \eqref{eq:chiBd} in parts, and treat each part separately with different methods:
\begin{equation}
    \sum_{\substack{u,v\in\Zd: \\ u\neq v \\ u-v\in r\Zd}}\tau_{\sZ,p}(u)\,\tau_{\sZ,p}(v)
    \le 2\sum_v \sum_{\substack{u\colon u\neq v \\ |u|\le|v|\\ u-v\in r\Zd}}\tau_{\sZ,p}(u)\,\tau_{\sZ,p}(v)\\
    = 2\Big((A)+(B)+(C)+(D)\Big),
\end{equation}
where
\begin{equation}
\begin{split}
    (A) = \sum_v \;\sum_{\substack{2r\le|u|\le|v|\\ u-v\in r\Zd}}\tau_{\sZ,p}(u)\,\tau_{\sZ,p}(v),\qquad
    &(B) = \sum_{|v|>MV^{1/6}\log V} \;\sum_{\substack{u\colon |u|\le2r\\ u-v\in r\Zd}}\tau_{\sZ,p}(u)\,\tau_{\sZ,p}(v)\\
    (C) = \sum_{2r<|v|\le MV^{1/6}\log V} \;\sum_{\substack{u\colon  |u|\le2r\\ u-v\in r\Zd}}\tau_{\sZ,p}(u)\,\tau_{\sZ,p}(v),\qquad
    &(D) = \sum_{|v|\le2r} \;\sum_{\substack{u\colon |u|\le2r\\ u-v\in r\Zd}}\tau_{\sZ,p}(u)\,\tau_{\sZ,p}(v)
\end{split}
\end{equation}
and $M$ is a (large) constant to be fixed later in the proof.
We proceed by showing that each of the four summands is bounded by a constant times $V^{-1/3}$, in that showing \eqref{eq:chiBd}.

Consider $(A)$ first. To this end, we prove for fixed $v\in\Zd$,
\begin{equation}\label{eq:chiBdProof1}
    \sum_{\substack{2r\le|u|\le|v|\\ u-v\in r\Zd}}\tau_{\sZ,p}(u)
    \le C_{\sss \tau}\,\frac{|v|^2}{V}.
\end{equation}
Indeed,
\begin{equation}\label{eq:chiBdProof2}
    \sum_{\substack{2r\le|u|\le|v|\\ u-v\in r\Zd}}\tau_{\sZ,p}(u)
    \le
    \sum_{\substack{2\le|u|\le\frac{|v|}{r}+1\\ u\in \Zd}}\tau_{p_c}\big(ru+(v\operatorname{mod} r)\big).
\end{equation}
By \refeq{powerbd}, this is bounded above by
\begin{equation}\label{eq:chiBdProof3}
    C_{\sss \tau}\,\sum_{2\le|u|\le\frac{|v|}{r}+1}\Big(r\big(|u|-1\big)+1\Big)^{-(d-2)}
    \le \frac{C_{\sss \tau}}{r^{d-2}}\,\sum_{1\le|u|\le\frac{|v|}{r}}|u|^{-(d-2)}.
\end{equation}
The discrete sum is dominated by the integral
\begin{equation}\label{eq:chiBdProof4}
    C_{\sss \tau}\,r^{-(d-2)}\int_{0\le|u|\le\frac{|v|}r}|u|^{-(d-2)}\,{\rm d}u
    \le C_{\sss \tau}\,C_\circ\,r^{-d}\,\frac{|v|^2}{2}
    \le C_{\sss \tau}\,C_\circ\,\frac{|v|^2}{V},
\end{equation}
as desired (with $C_\circ$ denoting the surface of the $(d-1)$-dimensional hypersphere).
Consequently, using \eqref{eq:chiBdProof1},
\begin{equation}\label{eq:chiBdProof5}
    (A)\le \frac{C_{\sss \tau}\,C_\circ}{V}\sum_v|v|^2\tau_{\sZ,p}(v)
    \le \frac{C_{\sss \tau}\,C_\circ}{V}\,\xi_2(p)^2\,\chi_\sZ(p)
    \le \frac{C_{\sss \tau}\,C_\circ\,\Cxit^2\,\Cg}{V}\,\big(p_c(\Zd)-p\big)^{-2}
\end{equation}
by the bounds in Theorem \ref{thm:ZdResults}.
Inserting $p=p_c(\Z^d)-{K}{\cn}^{-1}V^{-1/3}$ yields the desired upper bound $(A)\le C\,V^{-1/3}$.

For the bound on $(B)$ we start by calculating
\begin{equation}\label{eq:chiBdProof11}
    \sum_{u\colon |u|\le2r}\tau_{\sZ,p}(u)
    \le
    \sum_{u\colon |u|\le2r}\tau_{p_c(\Zd)}(u)
    \le
    \sum_{u\colon |u|\le2r}\frac{C_{\sss \tau}}{(|u|+1)^{d-2}}
    \le
    O(r^2).
\end{equation}
For the sum over $v$ we use the exponential bound of Theorem \ref{thm:ZdResults}:
From \refeq{defXi}--\refeq{xibd} and our choice of $p$ it follows that
$\tau_{\sZ,p}(v)\le\exp\big\{-C\,|v|\,V^{-1/6}\big\}$
for some constant $C>0$.
Consequently,
\begin{equation}\label{eq:chiBdProof12}
    \sum_{\substack{|v|>MV^{1/6}\log V\\ u-v\in r\Zd}}\tau_{\sZ,p}(v)
    \le
    \sum_{|v|>\frac M r V^{1/6}\log V}\tau_{\sZ,p}\big(rv+(u\operatorname{mod} r)\big)
    \le
    \sum_{|v|>\frac M r V^{1/6}\log V}\exp\big\{-r\big(|v|-1\big)CV^{-1/6}\big\}.
\end{equation}
This sum is dominated by the integral
\begin{equation}\label{eq:chiBdProof13}
    \int_{|v|>\frac M r V^{1/6}\log V}\exp\big\{-r\,|v|\,CV^{-1/6}\big\}
    \;\exp\big\{r\,C\,V^{-1/6}\big\}\,{\rm d}v,
\end{equation}
which can be shown by partial integration as being less or equal to
\begin{equation}\label{eq:chiBdProof14}
    \operatorname{const}(C,M,d)\,\frac{V^{d/6}}V\left(\log V\right)^d\exp\!\left\{-\frac M C \log V\right\}
    \;\exp\big\{r\,C\,V^{-1/6}\big\}.
\end{equation}
This expression equals
\begin{equation}\label{eq:chiBdProof15}
    \operatorname{const}(C,M,d)\,V^{d/6-1-M/C+C(1/d-1/6)}\,\left(\log V\right)^d.
\end{equation}
We now fix $M$ large enough such that the exponent of $V$ is less than $-(1/3+2/d)$.
This finally yields
\begin{equation}\label{eq:chiBdProof10}
    (B)
    \le
    \sum_{u\colon |u|\le2r} \;\sum_{\substack{|v|>MV^{1/6}\log V\\ u-v\in r\Zd}}\tau_{\sZ,p}(u)\,\tau_{\sZ,p}(v)
    \le
    \operatorname{const}(C,M,d)\,r^2\,o\big(V^{-(1/3+2/d)}\big)
    =
    o\big(V^{-1/3}\big).
\end{equation}

In order to bound $(C)$ we proceed similarly by bounding
\begin{equation}\label{eq:chiBdProof6}
\begin{split}
    (C)
    &\le C_{\sss \tau}^2\sum_{u\colon |u|<2r}\left(|u|+1\right)^{-(d-2)}
        \sum_{\substack{2r\le|v|\le MV^{1/6}\log V\\u-v\in r\Zd}}\left(|v|+1\right)^{-(d-2)}.
\end{split}
\end{equation}
A domination by integrals as in \eqref{eq:chiBdProof2}--\eqref{eq:chiBdProof4} allows for the upper bound
\begin{equation}\label{eq:chiBdProof7}
    C\; r^2 \; \frac{M^2\,V^{1/3}\,(\log V)^2}{V},
\end{equation}
and this is $o\big(V^{-1/3}\big)$ if $d>6$ for any $M>0$.

The final summand $(D)$ is bounded as in \eqref{eq:chiBdProof6} by
\begin{equation}\label{eq:chiBdProof8}
    C_{\sss \tau}^2\sum_{u\colon |u|<2r}\left(|u|+1\right)^{-(d-2)}
        \sum_{\substack{v\colon |v|\le 2r\\u-v\in r\Zd}}\left(|v|+1\right)^{-(d-2)}.
\end{equation}
The second sum can be bounded uniformly in $u$ by
\begin{equation}\label{eq:chiBdProof9}
    \sum_{\substack{v\colon |v|\le 2r\\u-v\in r\Zd}}\left(|v|+1\right)^{-(d-2)}
    \le
    (2r)^{-(d-2)}\;\#\{v\colon |v|\le 2r, \,u-v\in r\Zd\}
    \le (2r)^{-(d-2)}\;5^d,
\end{equation}
while the first sum is bounded by $C\,r^2$.
Together, this yields the upper bound $C\,r^{-(d-4)}$, and this is $o\big(V^{-1/3}\big)$ for $d>6$.

Finally, we have proved that $(A)\le C\,V^{-1/3}$, and that $(B)$, $(C)$, $(D)$ are of order $o\big(V^{-1/3}\big)$. This completes the proof of Lemma \ref{lemma:SumTauBound}.
\qed

\proof[Proof of Theorem \ref{thm-ImprovedBoundsCmax}.]
Assume that the conditions of Theorem \ref{thm-1} are satisfied.
Then by \cite[Corol.\ 4.1]{HeydeHofst07} there exists a constant $\Lambda>0$ such that, when $r\to\infty$,
\begin{equation}\label{eq:pcptub}
    p_c(\Zd)-p_c(\torus)\le\frac\Lambda\cn V^{-1/3}.
\end{equation}
It therefore suffices to prove a matching lower bound.

We take $p=p_c(\Z^d)-{K}{\cn}^{-1}V^{-1/3}$.
The following bound is proven in \cite{HeydeHofst07}:
\begin{equation}\label{eq:chiTbd}
    \chi_\sT(p)\quad\ge\quad\chi_\sZ(p)\,\Bigg(1-\bigg(\frac12+p\,\cn^2\,\chi_\sZ(p)\bigg)
    \sum_{\substack{u,v\in\Zd, u\neq v \\ u-v\in r\Zd}}\tau_{\sZ,p}(u)\,\tau_{\sZ,p}(v)\Bigg).
\end{equation}
Indeed, this bound is obtained by substituting \cite[(5.9)]{HeydeHofst07} and \cite[(5.13)]{HeydeHofst07} into \cite[(5.5)]{HeydeHofst07}.
Furthermore, by our choice of $p$ and \refeq{defCg},
$K^{-1}V^{1/3}\le\chi_\sZ(p)\le\Cg K^{-1}V^{1/3}$.
Together with \eqref{eq:chiBd},
\begin{equation}\label{eq:chiTbd2}
    \chi_\sT(p)\quad\ge\quad K^{-1}V^{1/3}\left(1-\left(1/2+p\,\cn^2\,K^{-1}\,\Cg\,V^{1/3}\right)
    \tilde C V^{-1/3}\right)
    \quad\ge\quad\tilde c_{\sss K}V^{1/3},
\end{equation}
where $\tilde c_{\sss K}$ is a small (though positive) constant.
Under the conditions of Theorem \ref{thm-1}, also the following bound holds by Borgs et al.\ \cite{BorgsChayeHofstSladeSpenc05a}:
For $q\ge0$,
    \eq\lbeq{SubcriticalPhase2}
        \chi_\sT\!\left(p_c(\torus)-\cn^{-1}q\right)\le\frac2q;
    \en
cf.\ the upper bound in \cite[(1.15)]{BorgsChayeHofstSladeSpenc05a}.
The upper bound \eqref{eq:pcptub} allows $K$ be so large that $p<p_c(\torus)$.
Consequently, the conjunction of \eqref{eq:chiTbd2} and \refeq{SubcriticalPhase2} obtains
    \eq
    \frac{2}{\cn(p_c(\torus)-p_c(\Zd)+KV^{-1/3})}
    \ge\chi_\sT(p)
    \ge \widetilde c_{\sss K} V^{1/3}.
    \en
This implies
    \eq\lbeq{sharpLowerBd2}
    p_c(\Zd)\ge p_c(\torus)+\left(K-\frac2{\widetilde c_{\sss K}\cn}\right) V^{-1/3},
    \en
as desired.
\qed

\vspace{.5em}
\noindent
The proof of Theorem \ref{thm-ImprovedBoundsCmax} concludes the proof of \refeq{Cmaxbd} in Theorem \ref{thm-1}, and it remains to prove \refeq{Cmaxbd2}.
\proof[Proof of \refeq{Cmaxbd2}]
The proof uses the exponential bound proven by Aizenman and Newman \cite[Proposition 5.1]{AizenNewma84}
that, for any $k\ge \chi_\sT(p)^2$,
\begin{equation}\label{eq:ExpProof1}
    \P_{\sT,p}\big(|\Ccal|\ge k\big)\le\left(\frac{{\operatorname{e}}}{k}\right)^{1/2}\exp\!\left\{-\frac{k}{2\,\chi_\sT(p)^2}\right\}.
\end{equation}
In order to apply this bound on the torus, we bound
\begin{equation}\label{eq:ExpProof2}
    \P_{\sT,p}\big(|\Cmax|\ge k\big)
    \quad\le\quad
    \frac1k\,\sum_{v\in\mathbbm{V}}\P_{\sT,p}\big(|\Cmax|\ge k,\,v\in\Cmax\big)
    \quad\le\quad
    \frac V k \;\P_{\sT,p}\big(|\Ccal|\ge k\big).
\end{equation}
Together with \eqref{eq:ExpProof1}, we obtain for $\omega>\chi_\sT(p)^2V^{-2/3}$,
\begin{equation}\label{eq:ExpProof4}
    \P_{\sT,p}\big(|\Cmax|\ge \omega V^{2/3}\big)
    \le\frac{{\rm e}^{1/2}}{\omega^{3/2}}\;\exp\!\left\{-\frac{\omega V^{2/3}}{2\,\chi_\sT(p)^2}\right\}.
\end{equation}
We now choose \ch{$p=p_c(\Zd)$ and use that $\chi_\sT(p_c(\Zd))\leq C_{\sss \chi}V^{1/3}$} to see that indeed, for $\omega>C_{\sss \chi}^2$, by \refeq{defCg},
\begin{equation}\label{eq:ExpProof5}
    \P_{\sT,p_c(\Zd)}\big(|\Cmax|\ge \omega V^{2/3}\big)
    \quad\le\quad
    \frac{{\rm e}^{1/2}}{\omega^{3/2}}\;
    \exp\!\left\{-\frac{\omega}{2\Cg^2}\right\}.
\end{equation}

\qed

\section{Proof of Theorem \ref{thm-2}}
\label{sec-thm-2}
\proof[Proof of \refeq{Cibd}.]
The upper bounds on $|\Ci|$ in Theorem \ref{thm-2} follow immediately
from the upper bound on $|\Cmax|$ in Theorem \ref{thm-1}. Thus, we only need
to establish the lower bound.

Recall the definition of $Z_{\sss \geq k}$ in \refeq{Zgeq-def}, and note that
    \eqn{
    \lbeq{ZExp}
    \Exp (Z_{\sss \geq k}) =V\, \P_{\sT, p}\big(|\C|\ge k\big).
    }
By construction, $|\Cmax|\geq k$ if and only
if $Z_{\geq k}\geq k$. We shall make essential use of properties of the
sequence of random variables $\{Z_{\sss \geq k}\}$ proved in
\cite{BorgsChayeHofstSladeSpenc05a}. Indeed,
\cite[\ch{Lemma 7.1}]{BorgsChayeHofstSladeSpenc05a}
states that, for all $p$ and all $k$, ${\rm Var}_p(Z_{\sss \geq k})\leq V\chi_{\sT}(p).$
When we take $p=p_c(\Zd)$, then,
by \eqref{eq:ClusterExpect} in Corollary \ref{cor:BCHSS-bounds}
above, there exists a constant $C_{\sss Z}$ such that $\chi_{\sT}(p_c(\Zd))\leq C_{\sss Z} V^{1/3}$.
Consequently,
    \eqn{\label{eq:VarUB}
    {\rm Var}_{p_c(\Zd)}(Z_{\sss \geq k})\leq C_{\sss Z} V^{4/3}
    }
uniformly in $k$.
Now, further, by \eqref{eq:DeltaEqualsTwo} in Corollary \ref{cor:BCHSS-bounds},
there exists $c_{\sss\C}>0$ such that
    \eqn{
    \P_{\sT, p_c(\Z^d)}\big(|\C|\ge k\big)\geq  \frac{2\,c_{\sss\C}}{\sqrt{k}}.
    }
Take $k=V^{2/3}/\omega$, for some $\omega\geq 1$ sufficiently large.
Together with the identity in \refeq{ZExp},
    \eqn{
    \lbeq{ZExp-rep}
    \expec_{p_c(\Zd)} (Z_{\sss \geq k}) \geq 2\,c_{\sss\C}\, \omega^{1/2} V^{2/3}.
    }
Thus, by the Chebychev inequality,
    \eqan{
    \prob_{p_c(\Zd)} \big(Z_{\sss \geq k}\leq c_{\sss\C} \omega^{1/2} V^{2/3}\big)
    &\leq \prob_{p_c(\Zd)} \Big(\big|Z_{\sss \geq k}-\expec_{p_c(\Zd)} (Z_{\sss \geq k})\big|
    \geq c_{\sss\C}\, \omega^{1/2} V^{2/3}\Big)\nn\\
    &\leq c_{\sss\C}^{-2}\, \omega^{-1}\, V^{-4/3}\,\VarPc(Z_{\sss \geq k})\;\leq\; \frac{C_{\sss Z}}{c_{\sss\C}^{2}\,\omega}.
    }
We take $\omega>0$ large. Then, the event
$Z_{\sss \geq k}> c_{\sss\C} \,\omega^{1/2} V^{2/3}$ holds with high probability.
On this event, there are two possibilities.
Either $|\Cmax|\geq c_{\sss\C}\,\omega^{1/2}V^{2/3}/i$,
or $|\Cmax|<c_{\sss\C}\,\omega^{1/2}V^{2/3}/i$, in which case there are at least $c_{\sss\C}\, \omega^{1/2} V^{2/3}/|\Cmax|
\geq i$ distinct clusters of size at least $k=\omega^{-1} V^{2/3}$. We conclude that
    \eqan{
    \Pbold_{\sT, p_c(\Z^d)}\big(|\Ci|\leq \omega^{-1}V^{2/3}\big)
    &\leq \prob_{p_c(\Zd)} \big(Z_{\sss \geq k}\leq c_{\sss\C}\, \omega^{1/2} V^{2/3}\big)
    +\prob_{p_c(\Zd)} \big(|\Cmax|\geq c_{\sss\C}\,\omega^{1/2}V^{2/3}/i\big)\nn\\
    &\leq \frac{C_{\sss Z}}{c_{\sss\C}^{2}\,\omega}+\frac{i \,\tilde b}{ c_{\sss\C}\omega},
    }
where $\tilde b$ is chosen appropriately from the exponential bound in \refeq{Cmaxbd2}.
This identifies $b_i$ as $b_i=i \tilde b /c_{\sss\C}+C_{\sss Z}/c_{\sss\C}^{2}$, and proves \refeq{Cibd}.
\qed
\vskip0.5cm

\noindent
We complete this section with the proof that any weak limit of $|\Cmax|V^{-2/3}$
is non-degenerate provided that a certain condition hold. 
Theorem \ref{thm-1} proves that the sequence $|\Cmax|V^{-2/3}$
is \emph{tight}, and, therefore, any subsequence of $|\Cmax|V^{-2/3}$
has a further subsequence that converges in distribution.

\begin{prop}[\ch{$|\Cmax|V^{-2/3}$ is not concentrated}]
\label{prop-limitsCmax}
Under the conditions of Theorem \ref{thm-1}, and assuming that there exists $\omega>6^{2/3}$ such that
	\eq
	\label{corr-cond}
	\liminf_{V\rightarrow \infty} V^{1/3} \Pbold_{\sT, p_c(\Z^d)}\Big(|\C|>\omega V^{2/3}\Big)>0,
	\en
the random sequence $|\Cmax|V^{-2/3}$ is non-concentrated.
\end{prop}

Mind that \eqref{eq:DeltaEqualsTwo} implies \eqref{corr-cond} for $\omega<b_{\sss \mathcal C}$ for a positive constant $b_{\sss \mathcal C}$ (which is the same as $b_1$ in \cite[Theorem 1.3]{BorgsChayeHofstSladeSpenc05a}). The actual value of $b_{\sss \mathcal C}$ depends of the position of $p_c(\Zd)$ within the critical window of $p_c(\torus)$. While Theorem \ref{thm-ImprovedBoundsCmax} guarantees that $p_c(\Zd)$ does lie within the critical window, it gives us no control on the precise position.

We believe that \eqref{corr-cond} is correct, in fact, even for \emph{all} $\omega>0$. 
For example, for the Erd\H{o}s-R\'enyi random graph model, which is the corresponding mean-field model, a corresponding statement is true for all $\omega>0$, cf.\ \cite[Lemma 2.2]{HofKagMul09}, where even a local limit version of \eqref{corr-cond} is proved.
However, we have not been able to show this for percolation on the high-dimensional torus.

In order to prove Proposition \ref{prop-limitsCmax}, we start by establishing a
\emph{lower bound} on the variance of $Z_{\sss \geq k}$. That is the content
of the following lemma:

\begin{lemma}[A lower bound on the variance of $Z_{\sss \geq k}$]
\label{lem-lower-bd-Var-Z}
For each $k\geq 1$,
    \eqn{
    \lbeq{lower-bd-Var-Z}
    \Var(Z_{\sss \geq k})\geq V\,\P_{\sT, p}\big(|\C|\ge k\big)\big[k-V\,\P_{\sT, p}\big(|\C|\ge k\big)\big].
    }
\end{lemma}

\proof We have that
    \eqn{
    \Var(Z_{\sss \geq k})=\sum_{u,v} \P_{\sT, p}\big(|\C(u)|\ge k,|\C(v)|\ge k\big)
    -\big[V\,\P_{\sT, p}\big(|\C|\ge k\big)\big]^2.
    }
Now, we trivially bound
    \eqn{
    \sum_{u,v} \P_{\sT, p}\big(|\C(u)|\ge k,|\C(v)|\ge k\big)
    \geq
    \sum_{u,v} \P_{\sT, p}\big(|\C(u)|\ge k,u\conn v\big)
    =V\,\expec[|\C|\indic{|\C|\ge k}]
    \geq V\,k\,\P_{\sT, p}\big(|\C|\ge k\big).
    }
Rearranging terms proves Lemma \ref{lem-lower-bd-Var-Z}.
\qed

\begin{lemma}[An upper bound on the third moment of $Z_{\sss \geq k}$]
\label{lem-upper-bd-third-mom-Z}
For each $k\geq 1$,
    \eqn{
    \lbeq{upper-bd-third-mom-Z}
    \expec_p[Z_{\sss \geq k}^3]\leq V\chi_\sT(p)^3+3\,\expec_p[Z_{\sss \geq k}]\,V\,\chi_\sT(p)+\expec_p[Z_{\sss \geq k}]^3.
    }
\end{lemma}

\proof We compute
    \eqan{
    \expec_p[Z_{\sss \geq k}^3]
    &=\sum_{u_1,u_2,u_3} \P_{\sT, p}\big(|\C(u_1)|\ge k,|\C(u_2)|\ge k, |\C(u_3)|\ge k\big)\nn\\
    &=\sum_{u_1,u_2,u_3} \P_{\sT, p}\big(|\C(u_1)|\ge k,u_1\conn u_2, u_3\big)\nn\\
    &\qquad +3\sum_{u_1,u_2,u_3} \P_{\sT, p}\big(|\C(u_1)|\ge k,u_1\conn u_2, |\C(u_3)|\ge k, u_1\nc u_3\big)\nn\\
    &\qquad+\sum_{u_1,u_2,u_3} \P_{\sT, p}\big(|\C(u_1)|\ge k,|\C(u_2)|\ge k, |\C(u_3)|\ge k, u_i\nc u_j\forall i\neq j\big)\nn\\
    &=(I)+3\,(II)+(III).
    }
We shall bound these terms one by one, starting with $(I)$,
    \eqan{
    (I)&\leq \sum_{u_1,u_2,u_3} \P_{\sT, p}\big(|\C(u_1)|\ge k,u_1\conn u_2, u_3\big)
    =V\,\expec_p[|\C|^2\indic{|\C|\geq k}]\leq V\,\expec_p[|\C|^2]\leq V\chi_\sT(p)^3,
    }
by the tree-graph inequality (see \cite{AizenNewma84}).
We proceed with $(II)$, for which we use the BK-inequality, to bound
    \eqan{
    (II)&\leq \sum_{u_1,u_2,u_3} \P_{\sT, p}(\{|\C(u_1)|\ge k, u_2\in \C(u_1)\}\circ\{|\C(u_3)|\ge k\}\big)\nnb
    &\leq \sum_{u_1,u_2,u_3} \P_{\sT, p}(|\C(u_1)|\ge k|, u_2\in \C(u_1))\,\P_{\sT, p}(\C(u_3)|\ge k) \nnb
    &=\;V\,\expec_p[|\C|\indic{|\C|\geq k}]\;\expec_p[Z_{\sss \geq k}]\leq \expec_p[Z_{\sss \geq k}]\,V\,\chi_\sT(p).
    }
We complete the proof by bounding $(III)$, for which we again use the BK-inequality, to obtain
    \eqan{
    (III)&\leq \sum_{u_1,u_2,u_3} \P_{\sT, p}(\{|\C(u_1)|\ge k\}\circ\{|\C(u_2)|\ge k\}\circ \{|\C(u_3)|\ge k\}\big)\nn\\
    &\leq \sum_{u_1,u_2,u_3} \P_{\sT, p}(|\C(u_1)|\ge k|)\,\P_{\sT, p}(\C(u_2)|\ge k)\,\P_{\sT, p}(|\C(u_3)|\ge k)
    =\expec_p[Z_{\sss \geq k}]^3.
    }
This completes the proof.
\qed
\vskip0.5cm

\noindent
Now we are ready to complete the proof of Proposition \ref{prop-limitsCmax}:
\proof[Proof of Proposition \ref{prop-limitsCmax}] By Theorem \ref{thm-1}, we know that
the sequence $|\Cmax|V^{-2/3}$ is tight, and so is $V^{2/3}/|\Cmax|$. Thus, there exists
a subsequence of $|\Cmax|V^{-2/3}$ that converges in distribution, and the weak limit,
which we shall denote by $X^*$, is strictly positive \ch{and finite} with probability 1.
Thus, we are left to prove that $X^*$ is non-degenerate. For this, we shall show that
there exists an $\omega>0$ such that $\prob(X^*>\omega)\in (0,1)$.

To prove this, we choose an $\omega$ that is not a discontinuity point of the distribution function of $X^*$ and note that
    \eqn{
    \prob(X^*>\omega)=\lim_{n\rightarrow \infty} \prob_{\sT, p_c(\Z^d)}(|\Cmax|V^{-2/3}_n>\omega),
    }
where the subsequence along which $|\Cmax|V^{-2/3}$ converges is denoted by $\{V_n\}_{n=1}^{\infty}$.
Now, using \refeq{Zgeq-def}, we have that
    \eqn{
    \prob_{\sT, p_c(\Z^d)}(|\Cmax|V^{-2/3}_n>\omega)
    =\prob_{\sT, p_c(\Z^d)}\big(Z_{\sss >\omega V_n^{2/3}}>\omega V_n^{2/3}\big).
    }
The probability $\prob_{\sT, p_c(\Z^d)}\big(Z_{\sss >\omega V^{2/3}}>\omega V^{2/3}\big)$
is monotone decreasing in $\omega$. By the Markov inequality and \eqref{eq:DeltaEqualsTwo},
for $\omega>6^{2/3}$ and uniformly in $V$,
    \eqn{
    \lbeq{Z*-tail-pre}
    \prob_{\sT, p_c(\Z^d)}\big(Z_{\sss >\omega V^{2/3}}>\omega V^{2/3}\big)
    \leq \omega^{-1} V^{-2/3}V\,\prob_{\sT, p_c(\Z^d)}\big(|\C|\ge \omega V^{2/3}\big)\leq \frac{6}{\omega^{3/2}}<1.
    }
In particular, the sequence $Z_{\sss >\omega V^{2/3}} V^{-2/3}$ is tight,
so we can extract a further subsequence \ch{$\{V_{n_l}\}_{l=1}^{\infty}$} so that also $Z_{\sss >\omega V^{2/3}} V^{-2/3}$
converges in distribution, say to $Z^*_{\omega}$. Then, \refeq{Z*-tail-pre} implies that
    \eqan{
    \lbeq{Z*-tail}
    \prob(Z^*_{\omega}=0)=1-\prob(Z^*_{\omega}>0)
    &=1-\lim_{l\to\infty} \prob_{\sT, p_c(\Z^d)}\big(Z_{\sss >\omega V^{2/3}_{n_l}}>0\big)\nn\\
    &=1-\lim_{l\to\infty} \prob_{\sT, p_c(\Z^d)}\big(Z_{\sss >\omega V^{2/3}_{n_l}}>\omega V^{2/3}_{n_l}\big)>0.
    }

Further, by Lemma \ref{lem-lower-bd-Var-Z},
    \eqan{
    \VarPc(Z_{\sss >\omega V^{2/3}} V^{-2/3})
    &\geq V^{-1/3}\P_{\sT, p_c(\Zd)}\big(|\C|>\omega V^{2/3})\big[\omega V^{2/3}-
    V\,\P_{\sT, p_c(\Zd)}(|\C|>\omega V^{2/3})\big]\nn\\  \label{eq:lowerBdVomega}
    &\geq V^{1/3}\P_{\sT, p_c(\Zd)}(|\C|>\omega V^{2/3})\big[\omega-6\omega^{-1/2}\big].
    }
We now apply \eqref{corr-cond} for suitable $\omega>6^{2/3}$ and the upper bound in \eqref{eq:DeltaEqualsTwo} to conclude that \eqref{eq:lowerBdVomega} is uniformly positive. 

Since there is also an upper bound on $\VarPc(Z_{\sss >\omega V^{2/3}} V^{-2/3})$ (this follows from \eqref{eq:VarUB}), it is possible to take a further subsequence $\{V_{n_{l_k}}\}_{k=1}^\infty$ for which
$\VarPc(Z_{\sss >\omega V^{2/3}} V^{-2/3})$ converges to $\sigma^2(\omega)>0$.
Since, by Lemma \ref{lem-upper-bd-third-mom-Z},
the third moment of $Z_{\sss >\omega V^{2/3}} V^{-2/3}$ is bounded, the random variable
$(Z_{\sss >\omega V^{2/3}} V^{-2/3})^2$ in uniformly integrable, and, thus,
along the subsequence for which $Z_{\sss >\omega V^{2/3}} V^{-2/3}$ weakly converges
and $\VarPc(Z_{\sss >\omega V^{2/3}} V^{-2/3})$ converges in distribution to $Z^*_{\omega}$,
we have
    \eqn{
    {\rm Var}(Z^*_{\omega})=\lim_{k\to\infty} \VarPc(Z_{\sss >\omega V^{2/3}_{n_{l_k}}} V_{n_{l_k}}^{-2/3})=\sigma^2(\omega)>0.
    }
Since ${\rm Var}(Z^*_{\omega})>0$, we must have that $\prob(Z^*_{\omega}=0)<1$.
Thus, by \refeq{Z*-tail} and the above, we obtain that $\prob(Z^*_{\omega}=0)\in (0,1),$ so that
    \eqan{
    \prob(X^*>\omega)&=\lim_{n\rightarrow \infty} \prob_{\sT, p_c(\Z^d)}\big(|\Cmax|V^{-2/3}_n>\omega\big)
    =\lim_{k\rightarrow \infty} \prob_{\sT, p_c(\Z^d)}\big(Z_{\sss >\omega V^{2/3}_{n_{l_k}}} V_{n_{l_k}}^{-2/3}>0\big)\nn\\
    &=\prob(Z^*_{\omega}>0)\in (0,1).
    }
This proves Proposition \ref{prop-limitsCmax}.
\qed

\section{Diameter and mixing time}
\label{sec-thm-3}
Let $d_{\Ccal}$ denote the graph metric (or \emph{intrinsic} metric) on the percolation cluster $\Ccal$.
\begin{theorem}[Nachmias--Peres \cite{NachmPeres08}]
\label{thm-NP}
Consider bond percolation on the graph $\gr$ with vertex set $\V$, $V=|\V|<\infty$, with percolation parameter $p\in(0,1)$.
Assume that for all subgraphs $\gr'\subset \gr$ with vertex set $\V'$,
\begin{itemize}
    \item[(a)]
    $\E_{\gr',p}\left|\edges\big(\{u\in\Ccal(v)\colon d_{\Ccal(v)}(v,u)\le k\}\big)\right|\le d_1 k,$ \qquad $v\in \V'$;
    \item[(b)]
    $\P_{\gr',p}\left(\exists u\in\Ccal(v)\colon d_{\Ccal(v)}(v,u)= k \right)\le d_2/k,$ \qquad $v\in \V'$,
\end{itemize}
where $\edges(\Ccal)$ denotes the number of open edges with both endpoints in $\Ccal$.
If for some cluster $\Ccal$
    \eq
    \label{eqCmaxbd}
    \Pbold_{\gr, p}\Big(\omega^{-1}V^{2/3}
    \leq |\Ccal|\Big)
    \geq 1-\frac{b}{\omega},
    \en
then there exists $c>0$ such that for all $\omega\ge1$,
\begin{eqnarray}\label{eq:Diam}
    \Pbold_{\gr, p}\Big(\omega^{-1}V^{1/3}
    \leq \diam(\Ccal)\leq \omega V^{1/3}\Big)
    &\geq& 1-\frac{c}{\omega^{1/3}},   \\
    \label{eq:TmixUB}
    \Pbold_{\gr, p}\Big(\Tmix(\Ccal)> \omega V\Big)
    &\leq& \frac{c}{\omega^{1/6}},  \\
    \label{eq:TmixLB}
    \Pbold_{\gr, p}\Big(\omega^{-1}V> \Tmix(\Ccal)\Big)
    &\leq& \frac{c}{\omega^{1/34}}.
\end{eqnarray}
\end{theorem}
We apply the theorem for $\gr=\torus$ and $p=p_c(\Zd)$.
Theorem \ref{thm-2} implies that \eqref{eqCmaxbd} holds for the $i$th largest cluster $\Ccal=\Ci$, $i\in\N$.
Hence Corollary \ref{thm-3} follows from Theorems \ref{thm-2} and \ref{thm-NP} once we have verified conditions (a) and (b) in the above theorem.
In fact, \eqref{eq:TmixUB} is a slight improvement over \eqref{eq:Tmix}.

Before proceeding with the verification, we shall comment on how to obtain Theorem \ref{thm-NP} from the work of Nachmias and Peres \cite{NachmPeres08}.
Indeed, Theorem \ref{thm-NP} is very much in the spirit of \cite[Theorem 2.1]{NachmPeres08}, though the $O$-notation there depends on $\beta$.
The bound \eqref{eq:Diam} is nevertheless straightforward from \cite[proof of Theorem 2.1(a)]{NachmPeres08} and \eqref{eqCmaxbd}.
For \eqref{eq:TmixUB} we use \eqref{eq:Diam} together with the bound
    $\Tmix(\graph)\le8\,|\edges|\,\diam(\graph)$,
valid for any finite (random or deterministic) graph $\graph$ with edge set $\edges$, cf.\ \cite[Corollary 4.2]{NachmPeres08}.

Furthermore, subject to conditions (a) and (b) of Theorem \ref{thm-NP},
there exist constants $C_1$, $C_2>0$ such that for any $\beta>0$, $D>0$,
\begin{equation}\label{eq:LB-Tmix}
    \P_{\gr,p}\!\left(\exists v\in\V\colon|\Ccal(v)|>\beta V^{2/3},\Tmix(\Ccal(v))< \frac{\beta^{21}}{1000\,D^{13}}V\right)
    \le D^{-1}\left(C_1+C_2\beta^3D^{-2}\right);
\end{equation}
which is obtained by combining \cite[(5.4)]{NachmPeres08} with the display thereafter.
From this we can deduce \eqref{eq:TmixLB} by choosing $D=1000^{-1/13}\omega$ and $\beta=\omega^{-1/34}$.

We complete the proof of Corollary \ref{thm-3} by verifying
that the conditions in Theorem \ref{thm-NP}(a) and (b) indeed
hold for critical percolation on the high-dimensional torus:

\paragraph{Verification of Theorem \ref{thm-NP}(a).}
The cluster $\C(v)$ is a subgraph of the torus with degree $\Omega$, therefore we can replace the number of edges on the left hand side by the number of vertices (and accommodate the factor $\Omega$ in the constant $d_1$).
In \cite[Proposition 2.1]{HeydeHofst07}, a coupling between
the cluster of $v$ in the torus and the cluster of $v$ in $\Z^d$
was presented, which proves that $\Ccal(v)$ can be obtained by
identifying points which agree modulo $r$ in a subset of the cluster
of $v$ in $\Z^d$. A careful inspection of this construction
shows that this coupling is such that it \emph{preserves}
graph distances. Since $\left|\{u\in\Ccal(v)\colon d_{\Ccal(v)}(v,u)\le k\}\right|$ is monotone in the
number of edges of the underlying graph, the result in Theorem \ref{thm-NP}(a)
for the torus follows from the bound $\Exp\left|\{u\in\Ccal(v)\colon d_{\Ccal(v)}(v,u)\le k\}\right|\le d_1 k$
for critical percolation on $\Z^d$. This bound was proved in
\cite[Theorem 1.2(i)]{KozmaNachm09}.
\qed

\paragraph{Verification of Theorem \ref{thm-NP}(b).}
For percolation on $\Z^d$, this bound was proved in \cite[Theorem 1.2(ii)]{KozmaNachm09}.
However, the event $\left\{\exists u\in\Ccal(v)\colon d_{\Ccal(v)}(v,u)= k \right\}$ is not monotone,
and, therefore, this does not prove our claim. However, a close inspection of the
proof of \cite[Theorem 1.2(ii)]{KozmaNachm09} shows that it only relies
on the bound that
    \eqn{
    \lbeq{delta-2}
    \prob_{\sT,p_c(\Zd)}(|\Ccal(v)|\geq k)\leq C_1/k^{1/2}
    }
(see in particular, \cite[Section 3.2]{KozmaNachm09}).
The bound \refeq{delta-2} holds for $k\le b_1 V^{2/3}$ by \cite[(1.19)]{BorgsChayeHofstSladeSpenc05a} and Theorem \ref{thm-ImprovedBoundsCmax} (where $b_1$ is a certain positive constant appearing in \cite[(1.19)]{BorgsChayeHofstSladeSpenc05a}).
For $k> b_1 V^{2/3}$ we use instead \eqref{eq:ExpProof1}.
Alternatively, one obtains \refeq{delta-2} from the corresponding $\Zd$-bound (proven by Barsky--Aizenman \cite{BarskAizen91} and Hara--Slade \cite{HaraSlade90a}), together with the fact that $\Zd$-clusters stochastically dominate $\torus$-clusters by \cite[Prop.\ 2.1]{HeydeHofst07}.
This completes the verification of Theorem \ref{thm-NP}(b).
\qed

\paragraph{Acknowledgement.}
The work of RvdH was supported in part by the Netherlands Organisation for Scientific Research
(NWO). We thank Asaf Nachmias for enlightening discussions concerning the results and
methodology in \cite{KozmaNachm09} and \cite{NachmPeres08}.
MH is grateful to Institut Mittag-Leffler for the kind hospitality during his stay in February 2009,
and in particular to Jeff Steif for inspiring discussions.

\def\cprime{$'$}


\begin{thebibliography}{10}

\bibitem{Aizen97}
M.~Aizenman.
\newblock On the number of incipient spanning clusters.
\newblock {\em Nuclear Phys. B}, 485(3):551--582, 1997.

\bibitem{AizenBarsk87}
M.~Aizenman and D.~J. Barsky.
\newblock Sharpness of the phase transition in percolation models.
\newblock {\em Comm. Math. Phys.}, 108(3):489--526, 1987.

\bibitem{AizenNewma84}
M.~Aizenman and C.~M. Newman.
\newblock Tree graph inequalities and critical behavior in percolation models.
\newblock {\em J. Statist. Phys.}, 36(1-2):107--143, 1984.

\bibitem{Aldou97}
D.~Aldous.
\newblock Brownian excursions, critical random graphs and the multiplicative
  coalescent.
\newblock {\em Ann. Probab.}, 25(2):812--854, 1997.

\bibitem{BarskAizen91}
D.~J. Barsky and M.~Aizenman.
\newblock Percolation critical exponents under the triangle condition.
\newblock {\em Ann. Probab.}, 19(4):1520--1536, 1991.

\bibitem{BorgsChayeHofstSladeSpenc05a}
C.~Borgs, J.~T. Chayes, R.~van~der Hofstad, G.~Slade, and J.~Spencer.
\newblock Random subgraphs of finite graphs. {I}. {T}he scaling window under
  the triangle condition.
\newblock {\em Random Structures Algorithms}, 27(2):137--184, 2005.

\bibitem{BorgsChayeHofstSladeSpenc05b}
C.~Borgs, J.~T. Chayes, R.~van~der Hofstad, G.~Slade, and J.~Spencer.
\newblock Random subgraphs of finite graphs. {II}. {T}he lace expansion and the
  triangle condition.
\newblock {\em Ann. Probab.}, 33(5):1886--1944, 2005.

\bibitem{Grimm99}
G.~Grimmett.
\newblock {\em Percolation}, volume 321 of {\em Grundlehren der Mathematischen
  Wissenschaften [Fundamental Principles of Mathematical Sciences]}.
\newblock Springer-Verlag, Berlin, second edition, 1999.

\bibitem{Hara90}
T.~Hara.
\newblock Mean-field critical behaviour for correlation length for percolation
  in high dimensions.
\newblock {\em Probab. Theory Related Fields}, 86(3):337--385, 1990.

\bibitem{Hara08}
T.~Hara.
\newblock Decay of correlations in nearest-neighbour self-avoiding walk,
  percolation, lattice trees and animals.
\newblock {\em Ann. Probab.}, 36(2):530--593, 2008.

\bibitem{HaraHofstSlade03}
T.~Hara, R.~van~der Hofstad, and G.~Slade.
\newblock Critical two-point functions and the lace expansion for spread-out
  high-dimensional percolation and related models.
\newblock {\em Ann. Probab.}, 31(1):349--408, 2003.

\bibitem{HaraSlade90a}
T.~Hara and G.~Slade.
\newblock Mean-field critical behaviour for percolation in high dimensions.
\newblock {\em Comm. Math. Phys.}, 128(2):333--391, 1990.

\bibitem{HaraSlade00a}
T.~Hara and G.~Slade.
\newblock The scaling limit of the incipient infinite cluster in
  high-dimensional percolation. {I}. {C}ritical exponents.
\newblock {\em J. Statist. Phys.}, 99(5-6):1075--1168, 2000.

\bibitem{HaraSlade00b}
T.~Hara and G.~Slade.
\newblock The scaling limit of the incipient infinite cluster in
  high-dimensional percolation. {II}. {I}ntegrated super-{B}rownian excursion.
\newblock {\em J. Math. Phys.}, 41(3):1244--1293, 2000.

\bibitem{HeydeHofst07}
M.~Heydenreich and R.~van~der Hofstad.
\newblock Random graph asymptotics on high-dimensional tori.
\newblock {\em Comm. Math. Phys.}, 270(2):335--358, 2007.

\bibitem{HeyHof17}
M.~Heydenreich and R.~van~der Hofstad.
\newblock {\em Progress in high-dimensional percolation and random graphs}.
\newblock {\it CRM Short Courses Series}. Springer, Cham, 2017.

\bibitem{HofKagMul09}
R.~van~der Hofstad, W.~Kager, and T.~M{\"u}ller.
\newblock A local limit theorem for the critical random graph.
\newblock {\em Electron. Commun. Probab.}, {\bf 14}:122--131, 2009.

\bibitem{KozmaNachm09}
G.~Kozma and A.~Nachmias.
\newblock The {A}lexander-{O}rbach conjecture holds in high dimensions.
\newblock {\em Invent. Math.}, 178(3):635--654, 2009.

\bibitem{Mensh86}
M.~V. Menshikov.
\newblock Coincidence of critical points in percolation problems.
\newblock {\em Dokl. Akad. Nauk SSSR}, 288(6):1308--1311, 1986.

\bibitem{NachmPeres09b}
A.~Nachmias and Y.~Peres.
\newblock Critical percolation on random regular graphs.
\newblock Preprint ar{X}iv:0707.2839v2 [math.{PR}], 2007.
\newblock To appear in \emph{Random Structures and Algorithms}.

\bibitem{NachmPeres08}
A.~Nachmias and Y.~Peres.
\newblock Critical random graphs: diameter and mixing time.
\newblock {\em Ann. Probab.}, 36(4):1267--1286, 2008.

\bibitem{Slade06}
G.~Slade.
\newblock {\em The Lace Expansion and its Applications}, volume 1879 of {\em
  Lecture Notes in Mathematics}.
\newblock Springer-Verlag, Berlin, 2006.

\end{thebibliography}
\end{document}